%% file: Main_Document/main.tex
\documentclass[12pt, shortlabels]{article}

\usepackage[english]{babel}
\usepackage[utf8]{inputenc}
\usepackage{lmodern}
\usepackage[T1]{fontenc}
\usepackage[normalem]{ulem}
\usepackage{bm}
\usepackage{fancyhdr}
\usepackage{enumitem}
\usepackage{aliascnt}

\usepackage{amsmath}
\usepackage{amsthm}
\usepackage{amssymb}
\usepackage{mathtools}
\usepackage{relsize}

\usepackage[export]{adjustbox}
\usepackage{graphicx}
\usepackage{wrapfig}
\usepackage{pdfpages}

\newtheorem{definition}{Definition}[subsection]
\newtheorem{example}[definition]{Example}
\newtheorem{theorem}[definition]{Theorem}
\newtheorem{lemma}[definition]{Lemma}

\DeclarePairedDelimiter\floor{\lfloor}{\rfloor}

\usepackage[a4paper]{geometry}
\usepackage{titlesec}
\usepackage{etoolbox}

\titleformat*{\section}{\Large \bfseries}
\titleformat*{\subsection}{\large \bfseries}
\titleformat*{\subsubsection}{\normalsize \bfseries}

\usepackage[hidelinks]{hyperref}
\usepackage[figure]{hypcap}
\usepackage{xurl}
\usepackage{cleveref}

\hypersetup{
    colorlinks,
    citecolor=red,
    filecolor=black,
    linkcolor=black,
    urlcolor=black
}

\begin{document}
\setcounter{page}{1}

\newcommand{\papertitle}{On the Existence of Optimal Strategies in a Combinatorial Game}
\newcommand{\authorname}{Tim Rammenstein}
\pagestyle{fancy}
\fancyhf{}
\fancyhead[C]{
  \footnotesize
  \textsc{
  \ifodd\value{page}
    \papertitle
  \else
    \authorname
  \fi
}}
\fancyfoot[C]{\footnotesize \textsc{\thepage}}

\renewcommand{\headrulewidth}{0pt}

\input{Section_1-Introduction/Introduction}
\input{Section_2-Foundations_of_the_Game/Foundations_of_the_Game}

\input{Section_3-Optimal_Strategies/Optimal_Strategies}

\input{Section_4-Conclusion/Conclusion}

\fancypagestyle{plain}{%
  \fancyhf{}
  \fancyhead[C]{%
    \footnotesize\textsc{%
      \ifodd\value{page}
        \papertitle
      \else
        \authorname
      \fi
    }%
  }
  \fancyfoot[C]{\footnotesize\textsc{\thepage}}
}

\thispagestyle{plain}

\thispagestyle{fancy}

\end{document}

%% file: Section_1-Introduction/Introduction.tex
\vspace*{2.5cm}

\begin{center}
    \LARGE  On the Existence of Optimal Strategies in a Combinatorial Game
    \vspace{0.8cm}

    \large Tim Rammenstein
    \vspace{0.3cm}
\end{center}

\begin{abstract} 
    \noindent
    We study a combinatorial game derived from a problem in the German National Mathematics Competition. In this game, two players take turns removing numbers from a finite set of natural numbers, aiming to satisfy a certain divisibility condition. We introduce a generalized version of the original game, which depends on two parameters: the size of the initial number set and a fixed divisor. For both players, we identify a broad range of game variants in which they can force a win. In particular, we show that for even-sized sets, the second player to move can always win, while for many odd-sized cases, the first player to move has a winning strategy. A web implementation of the game demonstrates some of our results in practice.
\end{abstract}

\section{Introduction} \label{Section 1}

In a problem from the 2017 \textit{Bundeswettbewerb Mathematik} (German National Mathematics Competition)$^{\text{\cite{num1}}}$, a game is described in which two players compete, and the task is to determine which player can force a win. During a project-based school week, the game was simplified so that it could be played quickly against others, while still allowing one player to have a guaranteed winning strategy.
\vspace{\lineskip}

The game can be further generalized by treating its two key parameters as variables. In this paper, we analyze the generalized variants of the game and rigorously justify the optimal strategies found in the course of analyzing the game. In doing so, statements are established that guarantee the existence of an optimal strategy for a wide range of game variants. Complementing the theoretical work, a website was created where variants of the game -- which I have named \textit{Zahlenschlacht} (this roughly translates to ``Number Battle``) -- can be played.
\vspace{\lineskip}

What makes the game Zahlenschlacht especially interesting are its simple gameplay and its short duration. Through finitary relations and number theory, particularly modular arithmetic, concepts of the game can be mathematically modeled and effectively studied.

%% file: Section_2-Foundations_of_the_Game/Foundations_of_the_Game.tex
\thispagestyle{plain}

\section{Foundations of the Game} \label{Section 2}
In this section, the game is introduced in its original form$^{\text{\cite{num1}}}$ as well as in a simplified variant.
After briefly situating the game within the domain of game theory, a generalization of the game is presented. Subsequently, fundamental game concepts are formalized mathematically, which is essential for the strategies and proofs in \hyperref[Section 3]{Section \ref*{Section 3}}.

\subsection{Rules of the Game}

The game Zahlenschlacht is a two-player game. Players take turns crossing out numbers displayed on a surface, such as a blackboard (hereafter referred to as the \emph{board}). The act of crossing out a number is referred to as a `move`.

\subsubsection{Original Game} \label{Original Game}

The numbers 1, 2, 3, \dots, 2017 are placed on the board. Player $A$ and Player $B$ take turns crossing out one number at a time until just two numbers remain. $A$ begins. If the sum of the two remaining numbers is divisble by 8, $A$ wins; otherwise, $B$ wins.$^{\text{\cite{num1}}}$
\vspace{\lineskip}

As shown in \cite{num2}, $A$ can force a win. For this, $A$ needs to play in accordance with the following strategy: $A$ crosses out, as the first move, a number that leaves a remainder of 1 when divided by 8, for example, the number 2017. Afterwards, for each number crossed out by $B$, $A$ crosses out a number such that the sum of the two numbers is divisible by 8.

\subsubsection{Game Variant} \label{Game Variant}

As part of a project-based school week, I searched for a simplified and easily playable variant of the original game in which one player can force a win. The following variant was developed in the process: The numbers 1, 2, 3, \dots, 15 are placed on the board. As in the original game, numbers are crossed out until just two remain. If the sum of the two remaining numbers is divisible by 7, $A$ wins. Otherwise, $B$ wins.

\begin{figure}[h]
    \centering
    \includegraphics[width=0.36\textwidth]{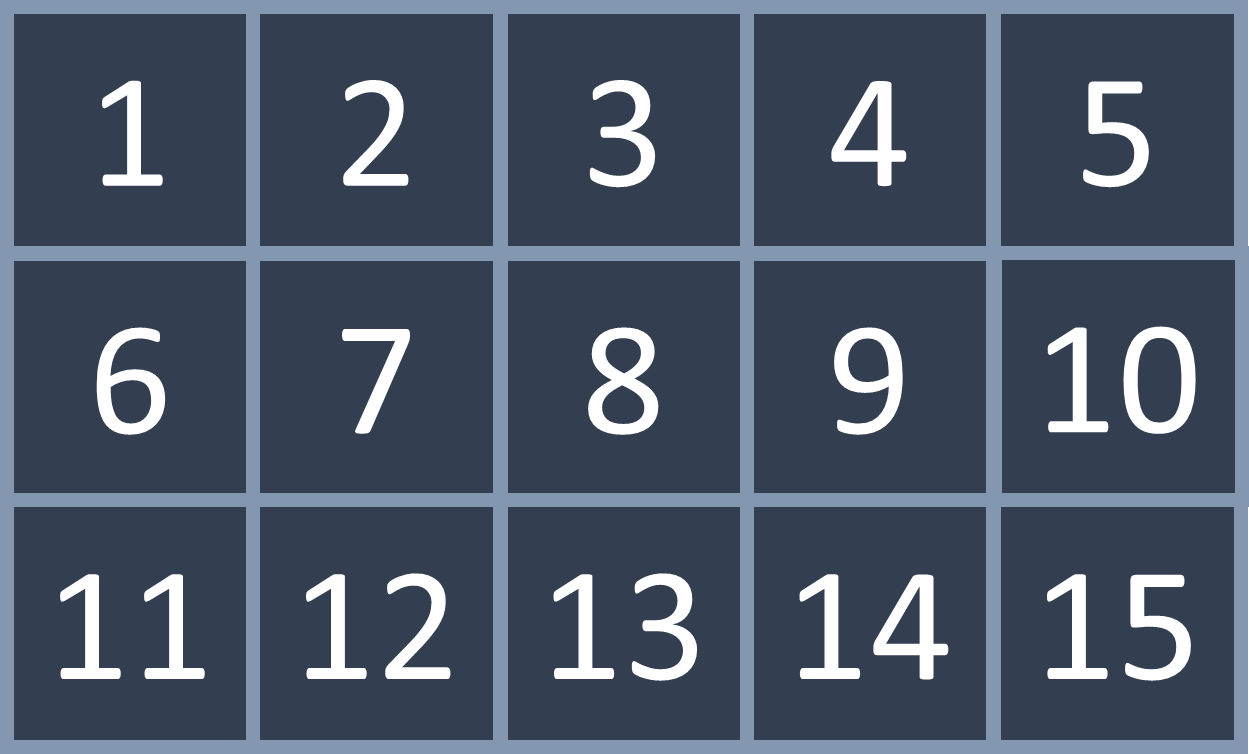}
    \caption{Visual representation of the board consisting of 15 numbers}
    \label{Fig:1}
\end{figure}

Here too, $A$ can force a win -- using a strategy similar to that in the original game: $A$ crosses out, as the first move, one of the numbers 1, 8 or 15. Then, after each move of $B$, $A$ crosses out a number such that the sum of both chosen numbers is divisible by 7. 
\vspace{\lineskip}

That $A$ can force a win using this strategy, is proven in \hyperref[Section 3]{Section \ref*{Section 3}}.

\subsection{Game-theoretic Classification}

In contrast to cooperative games, where players can form binding agreements to maximize a common benefit, players in non-cooperative games act independently and pursue their own interests. Zahlenschlacht belongs to the latter category.
\vspace{\lineskip}

The game is played as follows: Two players transform an initial position into a final position through a finite number of alternating moves, at which point a unique winner is determined (a draw is excluded, as the sum of the last two remaining numbers is either divisible by a fixed number or  not). Accordingly, Zahlenschlacht falls into the category of \emph{combinatorial games}$^{\text{\cite{num3}}}$ -- deterministic two-player games with alternating, well-defined moves and perfect information. Furthermore, it is an impartial game, as Player $A$ and Player $B$ can make the same moves from the same game position.

\subsection{Generalization of the Game} \label{2.3}

The difference between the game variant in \ref{Game Variant} and the original game in \ref{Original Game} lies only in the number of elements on the board and the number that decides the winner based on whether the sum of the last two numbers is divisible by it.
\vspace{\lineskip}

Let $n \in \mathbb{N}$ be the number of elements on the board. Clearly, it is reasonable to exclude $n = 1$, since in this case no two numbers could remain at the end. Likewise, $n = 3$ is also not to be considered, as the game would consist of only a single move and thus would not qualify as a two-player game.
\vspace{\lineskip}

In \hyperref[3.2]{Section \ref*{3.2}}, it is shown that for any even $n$, Player $B$ can force a win. Since this result holds for every second natural number, the focus of this paper lies on the considerably more interesting variants with odd $n$.
\vspace{\lineskip}

Let $d \in \mathbb{N}$ be the number that decides the winner based on whether the sum of the last two numbers is divisible by it. We assume $d \geq 2$, since for $d = 1$, Player $A$ trivially wins.
\vspace{\lineskip}

The parameters $n$ and $d$ define a generalization of the game from \cite{num1}, enabling the formulation of all meaningful variants.

\begin{definition} \label{General Rules}
Let $n$, $d$ $\in \mathbb{N}$, with $n \geq 4$ and $d \geq 2$, be fixed. We denote by $\textit{Z}(n, d)$ the two-player game with the following rules.
\end{definition}

\begin{itemize}
   \item \textit{The board consists of the first $n$ natural numbers in ascending order, $1$, \dots, $n$.}
    \item \textit{Player $A$ and Player $B$ take turns crossing out numbers until only two numbers remain. Player $A$ begins.}
    \item \textit{If the sum of the two remaining numbers is divisible by $d$, Player $A$ wins. Otherwise, Player $B$ wins.}
\end{itemize}

Unless stated otherwise, the numbers $n$ and $d$ are defined according to \autoref{General Rules} throughout the remainder of this paper.

\begin{example}
\emph{\hyperref[Original Game]{Section \ref*{Original Game}}} discusses the game $\textit{Z}(2017, 8)$, \emph{\hyperref[Game Variant]{Section \ref*{Game Variant}}} discusses the game $\textit{Z}(15, 7)$. In both games, $A$ can force a win. 
\end{example}

The game Zahlenschlacht encompasses the class of games $\textit{Z}(n, d)$. Rather than speaking of separate games, we also say that different variants of Zahlenschlacht are considered through $\textit{Z}(n, d)$ for various values of $n$ and $d$.

\subsection{Mathematical Formalization}

We will use modular arithmetic to reduce the game $\textit{Z}(n, d)$ to 
a mathematical model. This model allows for an effective analysis of Zahlenschlacht to identify and justify optimal strategies.
\vspace{\lineskip}

At each game state of the variant $\textit{Z}(n, d)$ of Zahlenschlacht, the board consists of a set $M \subseteq \{1, \dots, n\}$ of numbers, with $\lvert M \rvert \geq 2$. Since the remainder of division can be uniquely determined, it follows that the map
$$ f_{d} \colon M \to \{0, \dots, d - 1\}, \hspace{0.3cm} a \mapsto a \bmod d,  $$
assigning to each number in $M$ its unique remainder modulo $d$, is well-defined. We wish to replace the numbers on the board by their image under the map $f_{d}$. This would allow for a simplified analysis of the game, as the board would consist of smaller numbers and we could here apply modular arithmetic.
\vspace{\lineskip}

Since $f_{d}(a) = a \bmod d$ for any number $a \in M$, $f_{d}$ maps two numbers to the same residue if and only if they are congruent modulo $d$. Because the sum of two numbers modulo $d$ is congruent to the sum of their residues, we can replace the numbers $1, \dots, n$ on the board by their residues modulo $d$ because the sum of any two of those numbers leaves the same remainder as the sum of their residues modulo $d$. Particularly, the sum of the numbers is divisible by $d$ if and only if the sum of their residues is divisible by $d$. 

This way, the size of the original numbers is no longer relevant because only their residue modulo $d$ matters.

\begin{example} \label{2.4.1}
We consider the game $\textit{Z}(15, 7)$ with $M = \{1, \dots, 15\}$ and the map
$$ f_{7} \colon \{1, \dots, 15\} \to \{0, \dots, 6\}, \hspace{0.3cm} a \mapsto a \bmod 7.$$ 
Each of the numbers $1, \dots, 15$ on the board is mapped to its residue modulo $7$.

\begin{figure}[h]
    \centering
    \includegraphics[width=0.26\textwidth]{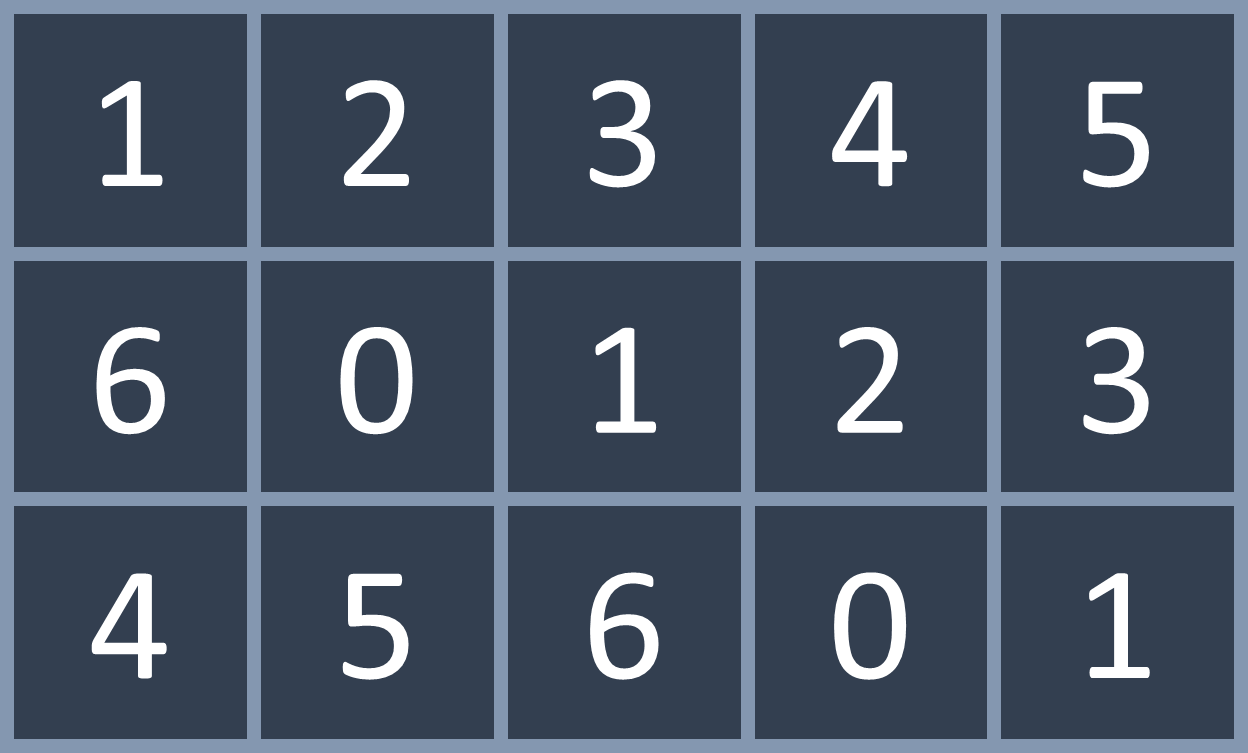}
    \caption{The residues of the numbers $1, \dots, 15$ modulo $7$}
    \label{Fig:2}
\end{figure}

The board containing the residues in \emph{\autoref{Fig:2}} is a mathematical model of \emph{\autoref{Fig:1}}, which reduces the numbers to their essential feature: their residues modulo $7$. 

According to \emph{\ref{Game Variant}}, the optimal strategy consists of crossing out one of the numbers $1, 8$ or $15$ in the first move. In this way, one number with residue $1$ modulo $7$ is crossed out. Without loss of generality, let $A$ cross out the number $1$.

\begin{figure}[h]
    \centering
    \includegraphics[width=0.26\textwidth]{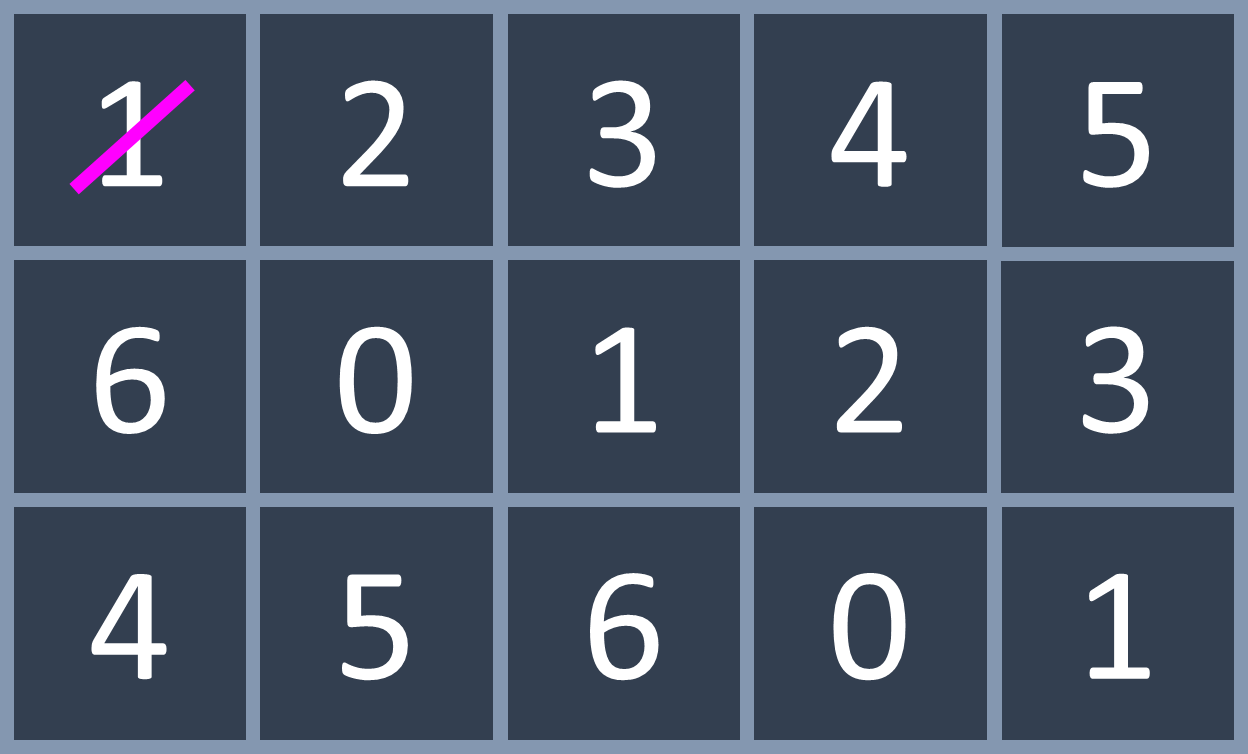}
    \caption{$A$ has crossed out the number 1 in the first move}
    \label{Fig:3}
\end{figure}

In \emph{\autoref{Fig:3}}, we can easily see that after this move of $A$, each of the remaining $14$ residues can be grouped into a pair with another residue such that the sum of both residues is divisible by $7$, and each residue on the board appears in exactly one of these pairs, which we will henceforth refer to as ``mod-$7$ pairs``. 
\vspace{\lineskip}

Since each residue is part of exactly one mod-$7$ pair, $B$ can only cross out a number from a complete mod-$7$ pair in the second move. Its corresponding partner can be crossed out by $A$. If, for example, $B$ crosses out a number with residue $2$ modulo $7$ $(the$ number $2$ or $9)$, A crosses out a number with residue $5$ modulo $7$ $(the$ number $5$ or $12)$ because $(2, 5)$ is a mod-$7$ pair.
\vspace{\lineskip}

As long as A follows this pattern for each move of B, it is ensured that only full mod-$7$ pairs are removed, so that the final two residues will also form such a pair. Thus their sum is divisible by $7$, leading to a win for $A$. 
\begin{figure}[h]
    \centering
    \includegraphics[width=0.26\textwidth]{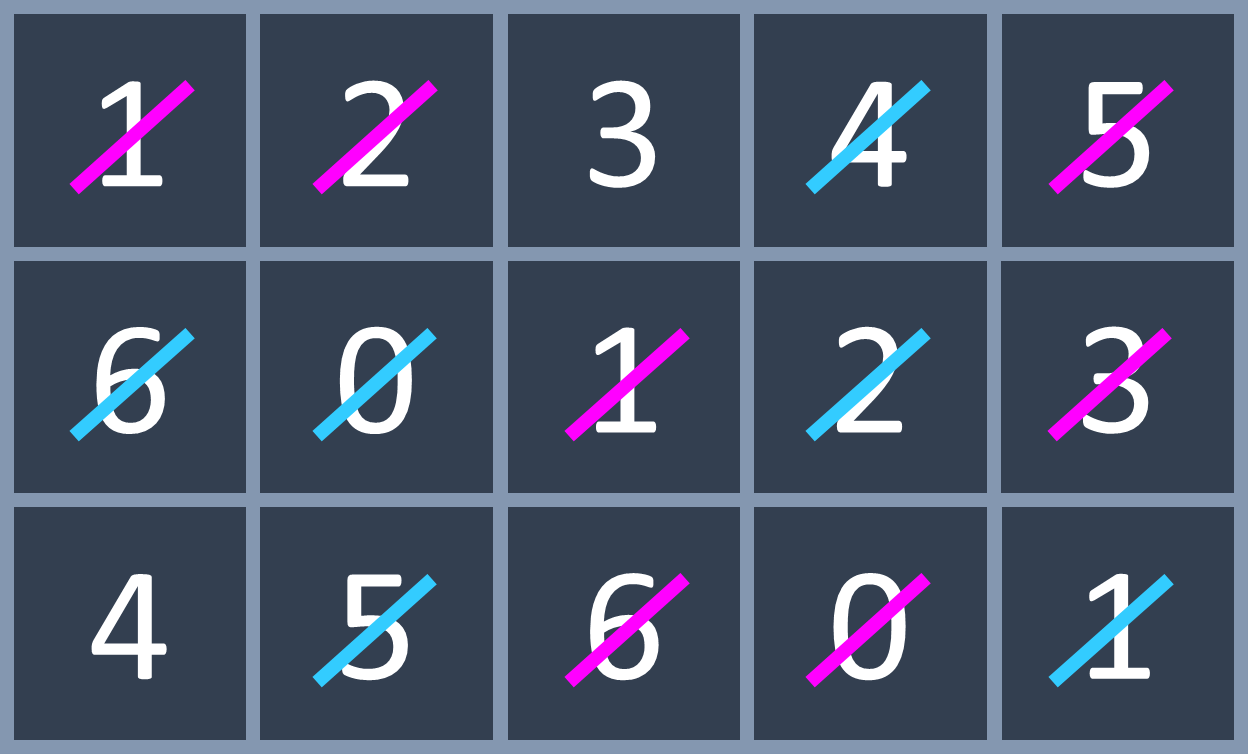}
    \caption{A possible course of play in which $A$ (pink) follows the strategy and wins against $B$ (turquoise)}
    \label{Fig:4}
\end{figure}
\end{example}

We generalize important concepts from \hyperref[2.4.1]{Example \ref*{2.4.1}} in \autoref{2.4.2}.

\begin{definition} \label{2.4.2}
For a given game state of the game $\textit{Z}(n, d)$, let $M \subseteq \{1, \dots, n\}$ with $\lvert M \rvert \geq 2$ be the set of numbers on the board.
\end{definition}

\begin{enumerate}
    \item \textit{The finitary relation $R_{d}$ on $\{1, \dots, n\}$ is the set of all pairs of distinct elements $x$, $y$ $\in \{1, \dots, n\}$, whose sum is divisible by $d$:
 $$ R_{d} := \{(x, y) \in \{1, \dots, n\}^{2}: x \neq y \land x + y \equiv 0 \mod d\}. $$}
    \item \textit{We call a number $a \in M$ \textbf{superfluous} if
 $$ \forall b \in M \setminus \{a\}: (a, b) \notin R_{d}. $$ 
 This means that $a$ is not related to any number in $M \setminus \{a\}$ by the relation $R_{d}$. }
    \item \textit{Let $r \in \{0, \dots, d - 1\}$ be an arbitrary residue modulo $d$. We denote by $a_{r}$ the number of numbers on the board that have residue $r$ modulo $d$:
    $$ a_{r} := \lvert \{x \in M: x \equiv r \mod d\} \rvert. $$}
    \item \textit{Let $\tilde{r}(n, d)$ denote the unique remainder of $n$ modulo $d$: $\tilde{r}(n, d) := n \bmod d$.}
\end{enumerate}

The relation $R_{d}$ is symmetric because it follows from $(x, y) \in R_{d}$ that $x \neq y$ and $y \neq x$ as well as $x + y = y + x \equiv 0 \mod d$, which means $(y, x) \in R_{d}$. So it is irrelevant whether we write $(x, y) \in R_{d}$ or $(y, x) \in R_{d}$.
\vspace{\lineskip}

In the model of the game $\textit{Z}(n, d)$ from \autoref{Fig:2} we can quickly count that, before the first move, that is, for $M = \{1, \dots, 15\}$, one has $a_1 = 3$ and $a_{r} = 2$ for all $r \in \{0, 2, 3, 4, 5, 6\}$. However, we can calculate the value for $a_{r}$, which is crucial for finding optimal strategies.

\begin{lemma} \label{2.4.3}
Before the first move, this means that $M = \{1, \dots, n\}$ is the set of numbers on the board of $\textit{Z}(n, d)$, we have for an arbitrary residue $r \in \{0, \dots, d - 1\}$ modulo d:
$$
a_{r} = \begin{cases}
  \floor*{\frac{n}{d}} + 1 &  \textit{for} \ 1 \leq r \leq \tilde{r}(n, d),  \\[0.2cm]
  \floor*{\frac{n}{d}} & \textit{for} \ r = 0 \lor r > \tilde{r}(n, d).
\end{cases}
$$
\end{lemma}

\begin{proof}
The greatest multiple of $d$ which is not greater than $n$, is $\floor*{\frac{n}{d}}d$. Each of the residues $0, \dots, d - 1$ modulo $d$ appears exactly $\floor*{\frac{n}{d}}$ times from 1 to $\floor*{\frac{n}{d}}d$. For all numbers lying between $\floor*{\frac{n}{d}}d$ and $n$, their respective residue appears one more time. Thus we have $a_{r} = \floor*{\frac{n}{d}} + 1$ for $1 \leq r \leq \tilde{r}(n, d)$. For the other residues $r$ that are smaller than 1 -- and thus satisfy $r = 0$ -- or greater than $\tilde{r}(n, d)$, we have $a_{r} = \floor*{\frac{n}{d}}$.
\end{proof}
\vspace{\lineskip}

\begin{example} \label{2.4.4}
We consider the game $\textit{Z}(15, 7)$ once more. We have $\tilde{r}(15, 7) = 15 \bmod 7 = 1$. Therefore, $1 \in \{0, \dots, 6\}$ is the only residue $r$ satisfying $1 \leq r \leq \tilde{r}(15, 7) = 1$. That is why, before the first move, we have $a_{1} = \floor*{\frac{15}{7}} + 1 = 2 + 1 = 3$ and $a_{r} = \floor*{\frac{15}{7}} = 2$ for all other residues $r \in \{0, 2, 3, 4, 5, 6\}$, which one can confirm in \emph{\autoref{Fig:2}} by counting.
\vspace{\lineskip}

$A$ wins if and only if, for the two remaining numbers $x, y$, we have $(x, y) \in R_{7}$. In \emph{\autoref{Fig:4}}, the number $3$ and the number $11$ which is replaced by its function value $f_{7}(11) = 4$ in the model, remain. Because we have $(3, 11) \in R_{7}$, $A$ wins.
\end{example}

%% file: Section_3-Optimal_Strategies/Optimal_Strategies.tex
\section{Optimal Strategies} \label{Section 3}

In this section, we analyze and justify optimal strategies for the game $\textit{Z}(n, d)$ that help one player to force a win. As mentioned in \ref{2.3}, $B$ can always win in games with an even $n$, which is proven in \ref{3.2}. That is why the statements in \ref{3.1} are only about games $\textit{Z}(n, d)$ with an odd $n$.

\subsection{Strategies for Player \textit{A}} \label{3.1}
At the heart of optimal strategies for Player $A$ lies a certain relation between the values of $a_{r}$. In this context, particular attention is paid to residues $\hat{r} \in \{0, \dots, d - 1\}$ that satisfy $\hat{r} + \hat{r} = 2\hat{r} \equiv 0 \mod d$. If twice $\hat{r}$ is divisible by $d$, then $\hat{r} = 0$, or, if $d$ is even, $\hat{r} = \frac{d}{2}$.

\begin{definition} \label{3.1.1}
In the game $\textit{Z}(n, d)$, let $M \subseteq \{1, \dots, n\}$ with $\vert M \vert \geq 2$ be the set of numbers on the board at a given game state.
\end{definition}

\begin{enumerate}
    \item \textit{The set $I_{d}$ contains exactly those residues $\hat{r} \in \{0, \dots, d - 1\}$ that are their own additive inverses modulo $d$: 
    $$ 
I_{d}  := \begin{cases}
  \left \{0, \frac{d}{2} \right\},  & \text{if} \ d \ \text{is even},  \\[0.2cm]
  \{0\}, & \text{if} \ d \ \text{is odd}.
\end{cases}
$$}
    \item \textit{We refer to a game state in which $a_{r} = a_{d - r}$ holds for all residues $r \in \{0, \dots, d - 1\} \setminus I_{d}$, and $a_{\hat{r}}$ is even for all $\hat{r} \in I_{d}$, as an \textbf{A-situation}.}
\end{enumerate}

The term ``$A$-situation`` and the notation using $a_{r}$ originate from \cite{num2} and are generalized in this paper. According to the following lemma, $A$ can secure a win if every number on the board can be paired with another number such that the sum of both numbers is divisible by 
$d$ and each number is part of exactly one such pair.

\begin{lemma} \label{3.1.2}
Let $n \geq 5$ be odd, $M \subseteq \{1, \dots, n\}$ with $\vert M \vert \geq 4$ be the set of numbers on the board at a given game state of the game $\textit{Z}(n, d)$, and let it be Player $B$’s turn. If an $A$-situation holds, Player $A$ can force a win. 
\end{lemma}

\begin{proof}
If $n$ is odd and $B$ is to move, the board contains an even amount of numbers. Let $B$ cross out an arbitrary number $a \in M$, and let $r := a \bmod d$ be the remainder of $a$ when divided by $d$. There are two cases for $r$:
\begin{enumerate}
    \item We have $r \in I_{d}$. Prior to $B$’s move, an $A$-situation was present, so $a_{r}$ was even; subsequently, $a_{r}$ is odd now. Then there is still at least one number with residue $r$ modulo $d$, which $A$ can cross out. Consequently, $a_{r}$ is even again. 
    \item We have $r \notin I_{d}$. Because an $A$-situation held prior to $B$’s move, we had $a_{r} = a_{d - r}$; now we have $0 \leq a_{r} = a_{d - r} - 1$, so $a_{d - r} \geq 1$ holds. Then $A$ can cross out a number with residue $d - r$ modulo $d$ so that $a_{r} = a_{d - r}$ holds again.
\end{enumerate}
In either case, the amount of numbers on the board has reduced by $2$ -- so it stays even -- and an \textit{A}-situation holds again. For all the other moves of $B$, $A$ can repeat this pattern until just two numbers remain, at which point an \textit{A}-situation will equally hold. It follows from the definition of an $A$-situation that in an $A$-situation with only two numbers, the sum of both numbers is congruent to $0 \mod d$ and thus divisible by $d$; so $A$ wins. 
\end{proof}  

\begin{example} \label{3.1.3}
We consider the games \emph{\ref{Original Game}} and \emph{\ref{Game Variant}}.
\end{example}

\begin{enumerate}
    \item \textit{In the game $\textit{Z}(2017, 8)$, we have $\tilde{r}(2017, 8) = 2017 \bmod 8 = 1$. According to \emph{\hyperref[2.4.3]{Lemma \ref*{2.4.3}}}, we thus have at the beginning $a_{1} = \floor*{\frac{2017}{8}} + 1 = 253$ and $a_{r} = \floor*{\frac{2017}{8}} = 252$ for all other residues $r \in\{0, \dots, 7\} \setminus \{1\}$. Since $8$ is even, $I_{8} = \{0, 4\}$ holds. If $A$ crosses out a number with residue $1$ modulo $8$ in the first move, we have afterwards $a_{r} = a_{8 - r} = 252$ for all $r \in \{0, \dots, 7\} \setminus I_{8}$. Moreover, $a_{0} = a_{4} = 252$ is even. Therefore, an $A$-situation is present prior to $B$’s next move, which means that $A$ can secure a win according to \emph{\hyperref[3.1.2]{Lemma \ref*{3.1.2}}}.}
    
    \item \label{2.} \textit{For $\textit{Z}(15, 7)$, we have $\tilde{r}(15, 7) = 15 \bmod 7 = 1$. Then $a_{1} = \floor{\frac{15}{7}} + 1 = 3$ and $a_{r} = \floor{\frac{15}{7}} = 2$ for all $r \in \{0, \dots, 6\} \setminus \{1\}$; particularly, $a_{0}$ is even. If $A$ removes a number with residue $1$ modulo $7$ in his next move -- that is, the number $1$, $8$ or $15$ -- an $A$-situation holds prior to $B$’s move, so $A$ can guarantee a win by \emph{\hyperref[3.1.2]{Lemma \ref*{3.1.2}}}.}
\end{enumerate}

With \hyperref[3.1.2]{Lemma \ref*{3.1.2}}, it is evident that $A$ can win the game with certainty whenever $A$ creates an $A$-situation in which it is $B$’s turn.
\vspace{\lineskip}

In the next step, we identify general game variants where $A$ can ensure that $B$ faces an $A$-situation. Here, $n$ is fixed, and $d$ depends on $n$. 

\begin{theorem} \label{3.1.4}
Let $n \geq 5$ be odd. Then $A$ can guarantee a win in the game $\textit{Z}(n, d)$ for all numbers $d$ with $d \in \left \{\frac{n - 1}{2}, \frac{n + 1}{2}, n, n + 1, n + 2 \right \}$.
\end{theorem}

\begin{proof}
We consider four cases for $d$. For each of these cases, we make use of \hyperref[3.1.2]{Lemma \ref*{3.1.2}}. 

\begin{enumerate}
    \item \label{1. of 3.1.4} Let $d = \frac{n - 1}{2}$. We have $\tilde{r}(n, d) = n \bmod \frac{n - 1}{2} = 1$, since $n = 2 \cdot \frac{n - 1}{2} + 1$. Additionally, $ \floor*{\frac{n}{\frac{n - 1}{2}}} = \floor*{\frac{2n}{n - 1}} = 2$ holds, because $2(n - 1) < 2n$ and $3(n - 1) > 2n$ for $n > 3$ and particularly for $n \geq 5$. Therefore, we have $a_{1} = 3$ before the first move, and $a_{r} = 2$ is even for all $r \in  \{0, \dots, d - 1 \} \setminus \left \{1 \right \}$. If $A$ crosses out a number with residue $1$ modulo $\frac{n - 1}{2}$, $B$ is faced with an $A$-situation before his move, so $A$ can secure a win.
    
    \item Let $d = n$. Then we have $a_{r} = 1$ for all $r \in \{0, \dots, d - 1\}$ in the beginning. If $A$ removes the number $n$ in the first move, $a_{0}$ is even. Since $d = n$ is odd, $I_{d} = \{0\}$ holds, so that $B$ then encounters an $A$-situation prior to his move. Consequently, $A$ can force a win. 

    \item Let $d = n + 1$. Then we have $\tilde{r}(n, d) = n \bmod (n + 1) = n$, and thus $a_{r} = \floor*{\frac{n}{n + 1}} + 1 = 0 + 1 = 1$ holds for all $r \in \{1, \dots, n\}$ as well as $a_{0} = 0$. Because $n$ is odd and $d = n + 1$ is thus even, $I_{d} = \{0, \frac{n + 1}{2}\}$ holds. Hence the only residue $i \in I_{d}$, for which $a_{i}$ is odd, is $i = \frac{n + 1}{2} \in I_{d}$. If $A$ removes the number $\frac{n + 1}{2}$ in the first move, $B$ is faced with an $A$-situation in the second move; with this, $A$ can force a win.

    \item Let $d = n + 2$. Since $\tilde{r}(n, d) = n \bmod (n + 2) = n$, we have $a_{r} = \floor*{\frac{n}{n + 2}} + 1 = 0 + 1 = 1$ for all $r \in \{1, \dots, n\}$ and $a_{0} = a_{n + 1} = 0$ in the beginning. To create an $A$-situation before it is $B$’s turn, $A$ crosses out the number $1$ in the first move, so that we then have: $a_{1} = 0 = a_{d - 1} = a_{(n + 2) - 1} = a_{n + 1}$. Because $d = n + 2$ is odd and $a_{0} = 0$ is already even, an $A$-situation is now present; thus $A$ can guarantee a win. 
\end{enumerate}

The statement for $d = \tfrac{n + 1}{2}$ follows immediately from the statement for $d = n + 1$, as $n + 1$ is, in particular, divisible by $\tfrac{n + 1}{2}$.
 \end{proof} 

\begin{example} \label{3.1.5}
Let $n = 15$. With \emph{\hyperref[3.1.4]{Theorem \ref*{3.1.4}}}, it follows that $A$ can secure a win in each of the games $\textit{Z}(15, 7)$, $\textit{Z}(15, 8)$, $\textit{Z}(15, 15)$, $\textit{Z}(15, 16)$ and $\textit{Z}(15, 17)$. 
\end{example}

We now define a strategy that can be used to transform some game situations into an $A$-situation in which $B$ is to move.

\begin{definition} \label{3.1.6}
Let $n \geq 5$ be odd, $d \geq 3$, and let $x, y, z \in \{0, \dots, d - 1\}$ be pairwise distinct. In the game $\textit{Z}(n, d)$, we refer to the following strategy for Player $A$ as the $\{x, y, z\}$\textbf{-strategy}:
\end{definition}

\begin{enumerate}
   \item \textit{$A$ crosses out in his first move a number $a \in M = \{1, \dots, n\}$ with arbitrary residue $r_{a} \in \{x, y, z\} \setminus I_{d}$ modulo $d$. Afterwards, $A$ responds to $B$’s moves in accordance with the strategy outlined in the proof of \emph{\hyperref[3.1.2]{Lemma \ref*{3.1.2}}}.}
    \item \textit{As soon as $B$ removes a number $b \in M \setminus \{a\}$ with residue $r_{b} \in \{x, y, z\} \setminus \{r_{a}\}$ for the first time, $A$ removes a number $c \in M \setminus \{a, b\}$ with residue $r_{c} \in \{x, y, z\} \setminus \{r_{a}, r_{b}\}$. Subsequently, $A$ again responds to $B$’s moves in accordance with the strategy outlined in the proof of \emph{\hyperref[3.1.2]{Lemma \ref*{3.1.2}}}.}
\end{enumerate}

The above-defined strategy cannot always be pursued. For example, in the game $\textit{Z}(15, 13)$, $B$ can prevent $A$ from pursuing the $\{0, 1, 2\}$-strategy by removing all numbers with a residue $r \in \{0, \dots, 12\} \setminus \{0, 1, 2\}$. The following lemma provides a sufficient condition for $A$ to pursue the $\{x, y, z\}$-strategy and thereby create an $A$-situation.

\begin{lemma} \label{3.1.7}
Let $n \geq 5$ be odd, $d \geq 3$, and let $A$ be to move in an arbitrary state of the game $\textit{Z}(n, d)$. Let $x, y, z \in \{0, \dots, d - 1\}$ be pairwise distinct residues that satisfy the following three conditions.
\end{lemma}
\begin{enumerate}[(1)]
    \item \label{Cond. 1} \textit{The values $a_{r}$ with $r \notin \{x, y, z\}$ and $(r, x), (r, y), (r, z) \notin R_{d}$ are precisely those that meet the criterion for an $A$-situation, and for all $s \in \{x, y, z\} \setminus I_{d}$, we have $a_{s} = a_{d - s} + 1$.}
    \item \label{Cond. 2} \textit{At least one of the three residues $x, y, z$ lies in $I_{d}$.}
    \item \label{Cond. 3} \textit{For all $i \in \{x, y, z\} \cap I_{d}$, we have $a_{i} \geq 3$.}
\end{enumerate}
\textit{Then $A$ can guarantee a win using the $\{x, y, z\}$-strategy.}

\begin{proof}
By \ref{Cond. 1}, only $a_{x}, a_{y}, a_{z}$ need to be adjusted for an $A$-situation to hold. Since $I_{d}$ has at most two elements, at least one of the residues $x, y, z$ does not lie in $I_{d}$, but at least one of them does, according to \ref{Cond. 2}. Without loss of generality, let $x \notin I_{d}$ and $z \in I_{d}$. By \ref{Cond. 1}, we then have $a_{x} = a_{d - x} + 1$. Now let $A$ remove a number with residue $x$ in the first move, so that $a_{x} = a_{d - x}$ meets the criterion for an $A$-situation. With that, only $a_{y}$ and $a_{z}$ are left to be adjusted. Since $z \in I_{d}$, it follows from \ref{Cond. 1} that $a_{z}$ is odd. It then suffices to cross out a number with residue $z$ to adjust $a_{z}$. Depending on whether $y \in I_{d}$ or $y \notin I_{d}$, the same principle applies to $y$ as to $z$ or $x$, respectively. Once $B$ removes a number with residue $y$ or $z$ for the first time, $A$ can respond by removing a number with residue $z$ or $y$, respectively, thereby creating an $A$-situation that secures a win for $A$.
\vspace{\lineskip}

It hence suffices to show that $B$ has to remove a number with residue $y$ or $z$ at some point in the game. If $B$ only crosses out numbers with a residue $r$, for which $a_{r}$ satisfies the criterion for an $A$-situation, $A$ can always remove a suited number to maintain the $A$-situation for these numbers. When all these numbers are removed, only numbers with a residue lying in $\{y, z, d - y\}$ remain. If $y \in I_{d}$, then $y$ and $d - y$ represent the same residue modulo $d$, and it follows from \ref{Cond. 3} that we have at least $3 + 3 = 6$ numbers on the board, so that $B$ is forced to cross out a number with residue $y$ or $z$. If $y \notin I_{d}$, $B$ can remove all numbers with residue $d - y$, and following the $\{x, y, z\}$-strategy, $A$ removes a number with residue $y$ in each turn. Because $a_{y} = a_{d - y} + 1$ initially holds in the case  $y \notin I_{d}$, one number with residue $y$ and at least three numbers with residue $z$ remain, hence at least four numbers in total. As $B$ then has at least one move remaining, $B$ must, at the latest at this point, remove a number with residue $y$ or $z$, enabling $A$ to force a win by pursuing the $\{x, y, z\}$-strategy.
\end{proof}

If two residues in $\{x, y, z\}$ lie in $I_{d}$, $A$ can also remove a number whose residue is contained in $I_{d}$ in the first move. Afterwards, the strategy stays the same.
\vspace{\lineskip}

The following theorem provides fixed values $d$ for which $A$ can secure a win in the game $\textit{Z}(n, d)$.

\begin{theorem} \label{3.1.8}
Let $n \geq 5$ be odd. 
\end{theorem}
\vspace{-0.6cm}

\begin{enumerate}[label=(\roman*)]
        \item \label{(i)} \textit{$A$ can guarantee a win in the game $\textit{Z}(n, d)$ for all $d \in \{2, 3\}$.}
        \item \label{(ii)} \textit{For $n \geq 11$, $A$ can guarantee a win in the game $\textit{Z}(n, d)$ for all $d \in \{2, 3, 4, 5, 6\}$.}
\end{enumerate}

\begin{proof}
We begin by showing \ref{(i)}. Let $d = 2$. Regardless of the moves made by $A$, when only three numbers remain on the board, the numbers are either all even, all odd, two even and one odd, or two odd and one even. In either case, $A$ can remove one of the numbers in such a way that the sum of the two remaining numbers is even and thus divisible by $d = 2$.
\vspace{\lineskip}

Now let $d = 3$. The numbers $n$, $n + 2$ and $n + 1$ are divisible by $3$ in the cases $n \bmod 3 = 0$, $n \bmod 3 = 1$, and $n \bmod 3 = 2$, respectively. Since by \hyperref[3.1.4]{Theorem \ref*{3.1.4}}, $A$ can secure a win in the games $\textit{Z}(n, n)$, $\textit{Z}(n, n + 2)$ and $\textit{Z}(n, n + 1)$, this also holds, in particular, for the game $\textit{Z}(n, 3)$.
\vspace{\baselineskip}

Now we show \ref{(ii)}. The statement for $d \in \{2, 3\}$ has been shown in the proof of \ref{(i)}. For the remaining values of $d$, we consider three cases.

\begin{enumerate}
    \item Let $d = 4$. Only the cases $n \bmod 4 \in \{1, 3\}$ need to be considered, as $n$ would be even otherwise. For $n \bmod 4 = 3$, $n + 1$ is divisible by $4$, and $A$ can secure a win by \hyperref[3.1.4]{Theorem \ref*{3.1.4}}. If $n \bmod 4 = 1$, then $\frac{n - 1}{2} = 2k$ for some $k \in \mathbb{N}$. For an even $k$, $2k$ is divisible by $4$, and $A$ can again force a win by \hyperref[3.1.4]{Theorem \ref*{3.1.4}}.
    For odd $k$, it follows from $n \bmod 4 = 1$ that, initially, $a_{1} = k + 1$ is even and $a_{0} = a_{2} = a_{3} = k$ is odd. Because the smallest odd $n \geq 11$, for which $n \bmod 4 = 1$ holds, is $n = 13$, and thus $\frac{n - 1}{2} = 2k \geq \frac{13 - 1}{2} = 6 = 2 \cdot 3$ holds, we have $a_{0} = a_{2} = k \geq 3$. It follows from \hyperref[3.1.7]{Lemma \ref*{3.1.7}} that $A$ can therefore pursue the $\{0, 1, 2\}$-strategy and secure a win.
    
    \item Let $d = 5$. For $n \bmod 5 = 1$, $\frac{n - 1}{2}$ is divisible by $5$, and for $n \bmod 5 \in \{0, 3, 4\}$, one of the numbers $n,$ $n + 2,$ $n + 1$ is divisible by $5$. The guaranteed win for $A$ follows again from \hyperref[3.1.4]{Theorem \ref*{3.1.4}}.
    \vspace{\lineskip}
    
   In the case $n \bmod 5 = 2$, we have $n = 5k + 2$ for some odd $k \in \mathbb{N}$, that is, $a_{1} = a_{2} = k + 1$ is even and $a_{0} = a_{3} = a_{4} = k$ is odd in the beginning. The smallest odd $n \geq 11$ with $n \bmod 5 = 2$ is $n = 17 = 5 \cdot 3 + 2$, so $a_{0} \geq 3$ holds. $A$ can then follow the $\{0, 1, 2\}$-strategy and thereby force a win again.
   
    \item Let $d = 6$. Because $n$ is odd, only the cases $n \bmod 6 \in \{1, 3, 5\}$ are to be considered. For $n \bmod 6 = 5$, $n + 1$ is divisible by $6$, and $A$ can secure a win by \hyperref[3.1.4]{Theorem \ref*{3.1.4}}. If $n \bmod 6 = 1$, there is some $k \in \mathbb{N}$ with $n = 6k + 1$. If $k$ is even, $a_{1} = k + 1$ is odd and $a_{i} = k$ is even for all $i \in \{0, 2, 3, 4, 5\}$ in the beginning. By removing a number with residue $1$ modulo $6$, $A$ can then guarantee a win by \hyperref[3.1.2]{Lemma \ref*{3.1.2}}. If $k$ is odd, the smallest such $n$ with $n \geq 11$ is the number $6 \cdot 3 + 1 = 19$. Thus we have $a_{0} = a_{3} = k \geq 3$. Then $A$ can make use of the $\{0, 1, 3\}$-strategy to secure a win.
    \vspace{\lineskip}
    
    If, on the other hand, $n \bmod 6 = 3$, then for an even $k \in \mathbb{N}$ satisfying $n = 6k + 3$, initially $a_{1} = a_{2} = a_{3} = k + 1$ is odd and $a_{0} = a_{4} = a_{5} = k$ is even. For an even $k$, the smallest such $n$ with $n \geq 11$ is $n = 6 \cdot 2 + 3 = 15$, so $a_{3} = k + 1 \geq 2 + 1 = 3$ holds. Then $A$ can force a win using the $\{1, 2, 3\}$-strategy. For an odd $k$ satisfying $n = 6k + 3 \geq 11$, we must have $k \geq 3$, and thus $a_{0} = k \geq 3$ holds. Then $A$ can pursue the $\{0, 1, 2\}$-strategy to secure a win.
\end{enumerate}
\vspace{-0.5cm}

\end{proof}

\begin{example} \label{3.1.9}
In the game $\textit{Z}(15, 6)$, initially $a_{1} = a_{2} = a_{3} = 3$ is odd and $a_{0} = a_{4} = a_{5} = 2$ is even. Using the $\{1, 2, 3\}$-strategy, $A$ can force a win. If, for example, $A$ removes the number $7$ in the first move, we then have $a_{1} = a_{5}$. If $B$ crosses out $2$ in the next move, $A$ should remove one of the numbers $3, 9, 15$, so that now not only $a_{2} = a_{4}$ holds, but also $a_{3} = 2$ is even. With that, an $A$-situation is present in which $B$ is to move; hence $A$ can force a win.
\end{example}

The following theorem provides formulas for $n$ that, for any fixed $d > 6$, generate infinitely many game variants in which $A$ can secure a win.

\begin{theorem} \label{3.1.10}
Let $d \geq 7$, and let $k \geq 1$. Then $A$ can guarantee a win in the game $\textit{Z}(n, d)$ for all $n$ with
\end{theorem}
\begin{enumerate}
    \item \label{Menge 1} $n \in \big \{kd - 1, \ (k + 1) \cdot d + 1, \ (k + 3) \cdot d - 3 \big \}$\textit{, if $d$ is even};
    \item \label{Menge 2} $n \in \big \{(2k - 1) \cdot d - 2, \ (2k - 1) \cdot d , \ 2kd - 1, \ 2kd + 1, \ (2k + 1) \cdot d + 2, \ (k + 1) \cdot 2d - 3 \big \},$ 
    
    \textit{if $d$ is odd}.
 \end{enumerate} 

\begin{proof}
If $d$ is even, the resulting numbers for $n$ are the sums or differences of a multiple of the even number $d$ and an odd number, so they are always odd. If $d$ is odd, $(2k \pm 1) \cdot d$ is odd too, as it is the product of two odd numbers. Note that addition or subtraction of $2$ does not change a number’s parity, so that $(2k - 1) \cdot d - 2, \ (2k - 1) \cdot d$ and $(2k + 1) \cdot d + 2$ are odd as well. The remaining values for $n$ that result from an odd $d$ are odd too, because they are the sum or difference of an even multiple of $d$ and an odd number.
\vspace{\lineskip}

Overall, all formulas yield only odd values for $n$, so that $A$ can apply the previously established strategies to the resulting game variants.

\begin{enumerate}
    \item Let $d$ be even. Since $d \geq 7$, we then have $d \geq 8$. By the formulas, we get $n \geq kd - 1 \geq d - 1 \geq 8 - 1 = 7 \geq 5$. This ensures that the game $\textit{Z}(n, d)$ is well-defined. If $n = kd - 1$, then $\tilde{r}(n, d) = d - 1$ holds, so that $a_{0} = \floor{\frac{kd - 1}{d}} = \floor{k - \frac{1}{d}} = k - 1$ and $a_{\frac{d}{2}} = a_{1} = \dots = a_{d - 1} = k$. For an even $k$, $A$ should remove a number with residue $0$. For an odd $k$, $A$ should remove a number with residue $\frac{d}{2}$. In either case, $A$ creates an $A$-situation, so $A$ can force a win.
    \vspace{\lineskip}

    If $n = (k + 1) \cdot d + 1$, then $a_{0} = a_{\frac{d}{2}} = a_{2} = \dots = a_{d - 1} = \floor{\frac{(k + 1) \cdot d + 1}{d}} = \floor{k + 1 + \frac{1}{d}} = k + 1$ and $a_{1} = k + 2$. If $k$ is even, we have $a_{0} = a_{\frac{d}{2}} = k + 1 \geq 2 + 1 = 3$, so that $A$ can follow the $\{0, 1, \frac{d}{2} \}$-strategy and force a win. Because, for an odd $k$, $a_{0} = a_{\frac{d}{2}} = k + 1$ is even, $A$ only needs to remove a number with residue $1$ in this case to create an $A$-situation and secure a win.
    \vspace{\lineskip}

    For $n = (k + 3) \cdot d - 3$, we have $\tilde{r}(n, d) = d - 3 \geq 8 - 3 = 5$, hence $2 < d - 3$, and, for $d > 6$, we have $\frac{d}{2} < d - 3$. Thus we have $a_{0} = a_{d - 2} = a_{d - 1} = \floor{\frac{(k + 3) \cdot d - 3}{d}} = \floor{k + 3 - \frac{3}{d}} = k + 2$ and $a_{\frac{d}{2}} = a_{1} = a_{2} = \dots = a_{d - 3} = k + 3$. If $k$ is even, then $a_\frac{d}{2} = k + 3 \geq 2 + 3 = 5$, so $A$ forces a win using the $\{1, 2, \frac{d}{2}\}$-strategy. If $k$ is odd, then $a_{0} = k + 2 \geq 1 + 2 = 3$, so that $A$ forces a win using the $\{0, 1, 2\}$-strategy.

    \item Let $d$ be odd. By the formulas, we get $n \geq (2 \cdot 1 - 1) \cdot d - 2 = d - 2 \geq 7 - 2 = 5$. Thus the game $\textit{Z}(n, d)$ is well-defined again. For $k = 1$, we get $n \in \{d - 2, \ d , \ 2d - 1, \ 2d + 1, \ 3d + 2, \ 4d - 3 \}$. Rewriting the first four values from this set in terms of $d$, we obtain the set $\{n + 2, n, \frac{n + 1}{2}, \frac{n - 1}{2} \}$. By \hyperref[3.1.4]{Theorem \ref*{3.1.4}}, we already know that $A$ can secure a win in game variants with such values for $d$.
    \vspace{\lineskip}

    We now consider $n = 3d + 2$. We have $a_{0} = a_{3} = \dots = a_{d - 1} = \floor{\frac{3d + 2}{d}} = \floor{3 + \frac{2}{d}} = 3$ and $a_{1} = a_{2} = 4$. Then $A$ can secure a win using the $\{0, 1, 2\}$-strategy.
    \vspace{\lineskip}

    Now let $n = 4d - 3$. We have $\tilde{r}(n, d) = d - 3$, so we get $a_{0} = a_{d - 2} = a_{d - 1} = \floor{\frac{4d - 3}{d}} = \floor{4 - \frac{3}{d}} = 3$ and $a_{1} = a_{2} = \dots = a_{d - 3} = 4$. Then $A$ can again force a win using the $\{0, 1, 2\}$-strategy.
    \vspace{\lineskip}

    Overall, $A$ can secure a win for all values of $n$ that are obtained by $k = 1$. Let $n_{0}$ be a value of $n$ computed by one of the formulas with $k = 1$, and let $a_{0}$ be the corresponding amount of numbers with residue 0 at the start of the game. For larger $k$, the values of $n$ generated by the same formula increase by an even multiple of $d$. Hence $a_{0}$ changes by an even number, preserving its parity. Furthermore, the value of $n \bmod d$ remains unchanged, so the partition into residues $r$ that satisfy $1 \leq r \leq \tilde{r}(n, d)$ and those that do not is the same across these values of $n$. Therefore, the initial configuration for all such $n$ is effectively identical. Since $A$ has a guaranteed winning strategy for $k = 1$, the same holds for all larger values of $k$. 
\end{enumerate}

\vspace{-0.5cm}

\end{proof}

\begin{example} \label{3.1.11}
Let $d = 42$. By the formulas, $A$ can secure a win in the game $\textit{Z}(n, 42)$ for all $n \in \{42k - 1, 42k + 43, 42k + 123\}$ with $k \geq 1$. Setting $k = 1$ yields the games $\textit{Z}(41, 42)$, $\textit{Z}(85, 42)$ and $\textit{Z}(165, 42)$; $k = 2$ yields $\textit{Z}(83, 42)$, $\textit{Z}(127, 42)$ and $\textit{Z}(207, 42)$.
\vspace{\lineskip}

For $d = 19$, we get $n \in \{38k - 21, 38k - 19, 38k - 1, 38k + 1, 38k + 21, 38k + 35\}$. Setting $k = 1$, we obtain $n \in \{17, 19, 37, 39, 59, 73\}$, and setting $k = 4$, we obtain $n \in \{131, 133, 151, 153, 173, 187\}$. For all these values of $n$, $A$ can secure a win in the game $\textit{Z}(n, 19)$.

If we only wish to compute values for $n$ with the same residue modulo $d$, we can simply choose one of the formulas for $n$. For example, the formula $n = 38k + 35$ for $d = 19$ yields exactly those values of $n$ that are congruent to $16$ modulo $19$, and for which $A$ can guarantee a win in the game $\textit{Z}(n, 19)$.
\end{example} 

\subsection{Strategies for Player \textit{B}} \label{3.2}

Player $B$ wins if the two final numbers in the game $\textit{Z}(n, d)$ are superfluous (see \autoref{2.4.2}). Hence the objective of $B$ is to make all numbers on the board superfluous, so that $A$ has no chance of winning.  While it is best for $A$ to remove superfluous numbers, with which $A$ can no longer reach a sum divisible by $d$ anyway, $B$ can ignore superfluous numbers and only needs to create more ones, until all numbers on the board are superfluous. Regardless of the numbers $A$ and $B$ remove at this point, $B$ wins.

\begin{theorem} \label{3.2.1}
Let $n \geq 5$ be odd. Then $B$ can guarantee a win in the game $\textit{Z}(n, d)$ for all numbers $d$ with $\frac{n + 3}{2} \leq d \leq n - 1$. 
\end{theorem}

\begin{proof}
Initially, we have $a_{0} = 1$ because $d < n$ and $2d \geq n + 3 > n$. Thus we have $a_{r} \in \{1, 2\}$ for all $r \in \{0, \dots, d - 1\}$. First, let $\frac{n + 3}{2} \leq d \leq n - 2$, which is only possible for $n \geq 7$. Then $\tilde{r}(n, d) = n \bmod d \geq 2$ holds. By \hyperref[2.4.3]{Lemma \ref*{2.4.3}}, it follows that $a_{1} = a_{2} = 2$ in the beginning. Moreover, we have $d - 1 > d - 2 \geq \frac{n + 3}{2} - 2 = \frac{n - 1}{2} > \frac{n - 3}{2} = n - \frac{n + 3}{2} \geq n \bmod d  = \tilde{r}(n, d)$. This means that $a_{d - 1} = a_{d - 2} = 1$ in the beginning. 

Since $d$ is superfluous due to $a_{0} = 1$, it is optimal for $A$ to remove either $d$ or one of the numbers with residue $1$ or $2$ in the first move, so that the condition $a_{i} = a_{d - i}$ is met for one index $i \in \{1, 2\}$. For any residue $s$ satisfying $a_{s} = 2$ and $a_{d - s} = 1$, an optimal response for $B$ is to remove the number $d - s$, thereby making both numbers with residue $s$ superfluous. In order to minimize the amount of superfluous numbers, $A$ should then always remove a superfluous number in his next move. However, for all such residues $s$, at least one new superfluous number is created. 

If $A$ removes $d$ in the first move, then $B$ can create at least one superfluous number for each of the residues $1$ and $2$. If instead $A$ initially removes a number with residue $1$ or $2$, then $B$ can create a superfluous number for the other residue in $\{1, 2\}$, and $d$ remains superfluous. In either case, there are at least two superfluous numbers. 

For any residue $r \notin I_{d}$ with $a_{r} = a_{d - r} = 1$, $B$ can remove either the number $r$ or $d - r$. Regardless of whether $A$ then removes a thereby newly created superfluous number or a number that was already superfluous, the total amount of superfluous numbers does not decrease. Since, as established, there are at least two superfluous numbers at this stage, there will also be two remaining at the end of the game. Their sum is certainly not divisible by $d$, and thus $B$ wins.
\vspace{\lineskip}

Now let $d = n - 1$. Then $d$ is even and we have $\tilde{r}(n, d) = n \bmod (n - 1) = 1$, so $a_{1} = 2$ and $a_{r} = 1$ for all $r \in \{0, 2, \dots, d - 1\}$ in the beginning. Hence there are initially two superfluous numbers, namely $\frac{d}{2} = \frac{n - 1}{2}$ and $d = n - 1$, and $A$ should remove one of them. If $B$ crosses out the number $d - 1$ afterwards, then the numbers $1$ and $n$ are superfluous. Now there are three superfluous numbers, one of which $A$ should cross out. After $A$’s move, there are at least two superfluous numbers left. For all remaining numbers, the corresponding residue $r$ satisfies $a_{r} = a_{d - r} = 1$. As before, $B$ can then make sure that the amount of superfluous numbers does not decrease, that is, its value is at least $2$. Therefore, $B$ can again force a win.
\end{proof}

\begin{example}
We consider the game $\textit{Z}(15, 9)$ with residues modulo $9$.

\begin{figure}[h]
    \centering
    \includegraphics[width=0.36\textwidth]{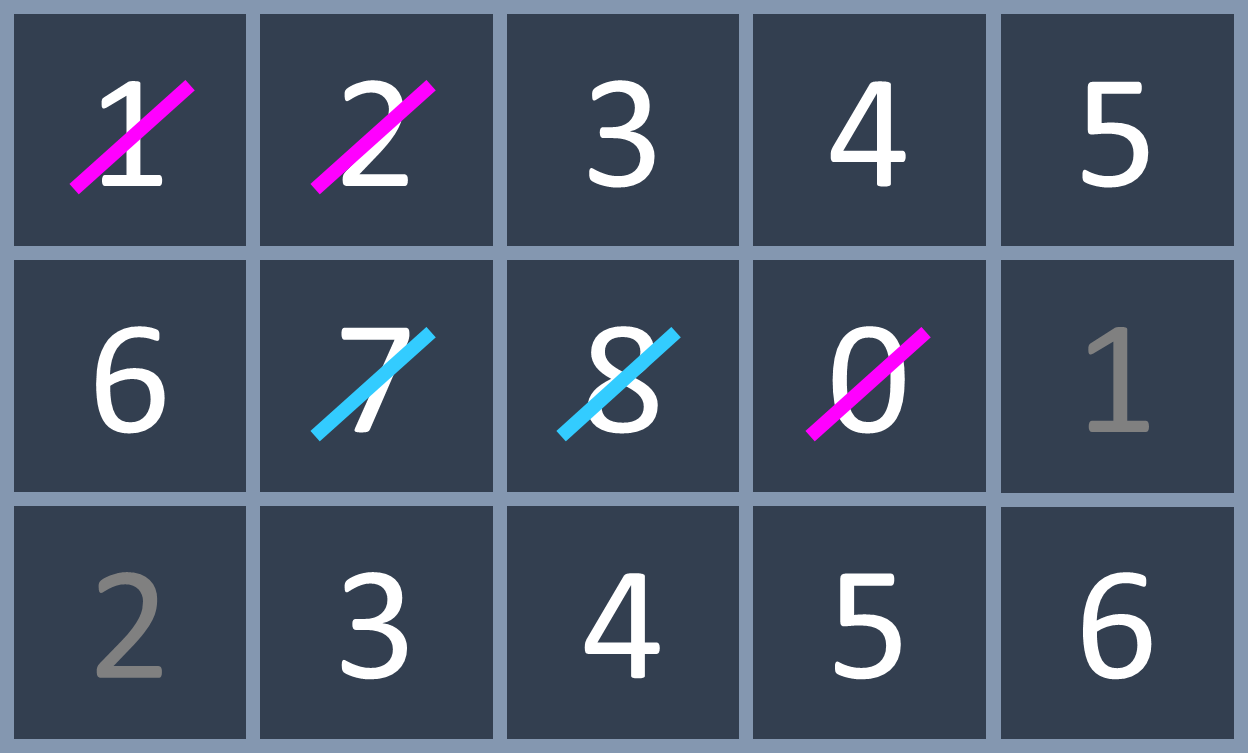}
    \caption{A course of play in which $A$ (pink) initially removes the number $9$ (residue $0$)}
    \label{Fig:5}
\end{figure}

In his first move, $B$ (turquoise) crosses out the number $8$, so that the numbers $1$ and $10$ -- both have residue $1$ -- become superfluous. Then $A$ can remove the number $1$, so that only $10$ is superfluous now. Next, $B$ proceeds analogously with the number $7$, and $A$ follows with the number $2$. Now the numbers $10$ and $11$, whose residues are marked in gray, are superfluous. For all $r \in \{3, 4, 5\}$, we have $a_{r} = a_{9 - r} = 2$. To secure a win, $B$ needs to make all numbers superfluous. In order to achieve this, for example, it suffices to remove all numbers with residues $3$ and $4$. 
\end{example}

By \hyperref[3.1.4]{Theorem \ref*{3.1.4}}, \hyperref[3.1.8]{Theorem \ref*{3.1.8}} and \hyperref[3.2.1]{Theorem \ref*{3.2.1}}, we know that for any fixed odd $n \geq 5$, $A$ can secure a win in the game $\textit{Z}(n, d)$ for $d \in \{2, 3, \frac{n - 1}{2}, \frac{n + 1}{2}, n, n + 1, n + 2\}$ and additionally for $d \in \{4, 5, 6\}$ if $n \geq 11$, while $B$ can secure a win for $d \in \{\frac{n + 3}{2}, \dots, n - 1\}$. However, for $d \in \{7, \dots, \frac{n - 3}{2}\}$, it cannot be determined in general which player can force a win: For example, choose $d = 7$. Here, $A$ can always win the game $\textit{Z}(29, 7)$ ($\frac{n - 1}{2} = \frac{29 - 1}{2} = 14$ is divisible by $d = 7$, so  $A$ can force a win by \hyperref[3.1.4]{Theorem \ref*{3.1.4}}), but $B$ can force a win in the game $\textit{Z}(17, 7)$ by removing all numbers with residues $4$, $5$ and $6$ and the number $7$ or $14$.
\vspace{\lineskip}

As stated in the next theorem, $B$ can force a win in the game variants with values of $d$ belonging to the last previously unexamined interval.

\begin{theorem} \label{3.2.3}
Let $n \geq 5$ be odd. Then $B$ can guarantee a win in the game $\textit{Z}(n, d)$ for all $d \geq n + 3$.    
\end{theorem}

\begin{proof}
For $d \geq 2n$, the statement is straightforward, because $n + (n - 1) = 2n - 1 < d$ is the largest possible sum of two numbers in the game, which is clearly not divisible by $d$. 

Now let $n + 3 \leq d \leq 2n - 1$, that is, $d = n + i$ for some $i \in \{3, \dots, n - 1\}$. Initially, the numbers $1, \dots, i - 1$ are always superfluous, since $d - 1 > \dots > d - (i - 1) = n + i - (i - 1) = n + 1 \notin \{1, \dots, n\}$. In other words, the smallest numbers that could form a sum divisible by $d$ together with the numbers $1, \dots, i - 1$ are never part of the numbers in the game. For $i \geq 4$, there are thus at least $4 - 1 = 3$ superfluous numbers from the very beginning. However, for $i = 3$, there are also three superfluous numbers, since $1$ and $2$ are included and $d = n + 3$ is even, so that $\frac{d}{2} = \frac{n + 3}{2} < n$ (because $n > 3$) and $\frac{d}{2} + d = \frac{3}{2}d = \frac{3}{2}(n + 3) > n$, implying that $a_{\frac{d}{2}} = 1$. It is best for $A$ if there are only three superfluous numbers in the beginning, so that after $A$’s first move, at most two superfluous numbers remain in the game. Apart from these two superfluous numbers, all other numbers have residues $r$ satisfying $a_{r} = a_{d - r} = 1$. Then $B$ can proceed as in the proof of \hyperref[3.2.1]{Theorem \ref*{3.2.1}}, ensuring that $A$ cannot reduce the amount of superfluous numbers. Hence the two final numbers will be superfluous, and $B$ wins.

Since there are at least three superfluous numbers at the beginning and $B$ can force a win in the worst-case scenario from his perspective -- where there are only three such numbers -- this holds, in particular, for all other cases. 
\end{proof} 

\begin{example} 
In the game $\textit{Z}(23, 26)$, the numbers $1$, $2$ and $13$ are superfluous. All other numbers lie in the set $\{1, \dots, 23\} \setminus \{1, 2, 13\}$ and are related to exactly one number by the relation $R_{26}$. In the first move, $A$ should remove one of the superfluous numbers, so that only two of them are left. For each remaining pair of numbers that are related by the relation $R_{26}$, $B$ can remove exactly one number in the corresponding pair. Consequently, all numbers become superfluous, so $B$ is guaranteed to win. The games $\textit{Z}(23, d)$ with $d \geq 27$ are similarly won by $B$.
\end{example}

As mentioned in \ref{2.3}, we now show that $B$ can always win when $n$ is even. This result clarifies why, up to this point, we have restricted our analysis to games with odd $n$.  

\begin{theorem} \label{3.2.5}
If $n$ is even, $B$ can guarantee a win in the game $\textit{Z}(n, d)$.
\end{theorem}

\begin{proof}
When $n$ is even, $B$ has the last move in the game. Then $B$ crosses out one of the three last remaining numbers, whose residues modulo $d$ we denote by $r_{1}, r_{2}, r_{3}$. If two of the three numbers have a sum not divisible by $d$, $B$ can remove the third one and win. Thus $A$ can only secure a win if the sum of any two of these numbers is divisible by $d$, that is, if $r_{1} + r_{2} \equiv r_{2} + r_{3} \equiv r_{1} + r_{3} \equiv 0 \mod d$ holds. The first congruence implies $r_{1} \equiv r_{3} \mod d$, and since $r_{1}$ and $r_{3}$ are (unique) residues modulo $d$, we have $r_{1} = r_{3}$. The second congruence implies $r_{2} \equiv r_{1} \mod d$, and by analogy, $r_{2} = r_{1}$. Overall, $r_{1} = r_{2} = r_{3}$ holds. Hence the sums imply that $2r_{i} \equiv 0 \mod d$ for $i \in \{1, 2, 3\}$, so $r_{1} = r_{2} = r_{3} \in I_{d}$. To secure a win, $B$ thus only needs to make sure that the final three numbers do not all have the same residue from the set $I_{d}$. Since $n$ is even and all but the two final numbers are removed, $A$ and $B$ each have $\frac{n - 2}{2}$ moves in total. Thus, before the last move, $B$ has $\frac{n - 2}{2} - 1 = \frac{n - 4}{2}$ moves during which $B$ must ensure that $a_{\hat{r}} < 3$ holds for all $\hat{r} \in I_{d}$. For $a_{0} < 3$ to hold, $B$ needs to remove at most $a_{0} - 2 = \floor*{\frac{n}{d}} - 2$ numbers with residue $0$.
\vspace{\lineskip}

Let $d$ be odd. Then $d \geq 3$, so $\floor*{\frac{n}{d}} - 2 \leq \floor*{\frac{n}{3}} - 2 \leq \frac{n}{3} - \frac{6}{3} = \frac{n - 6}{3} = \frac{2n - 12}{6} < \frac{3n - 12}{6} = \frac{n - 4}{2}$. Therefore, the $\frac{n - 4}{2}$ moves that $B$ has before the last move suffice to cross out enough numbers with residue $0$.
\vspace{\lineskip}

Let $d$ be even. Since $a_{\frac{d}{2}} \leq a_{0} + 1 = \floor*{\frac{n}{d}} + 1$, $B$ has to remove at most $a_{\frac{d}{2}} - 2 \leq \bigl( \floor*{\frac{n}{d}} + 1 \bigr) - 2 = \floor*{\frac{n}{d}} - 1$ numbers with residue $\frac{d}{2}$, for $a_{\frac{d}{2}} < 3$ to hold. In total, $B$ hence needs to remove at most $\big( \floor*{\frac{n}{d}} - 2 \big) + \big( \floor*{\frac{n}{d}} - 1 \big) = 2 \cdot \floor*{\frac{n}{d}} - 3$ numbers, for $a_{\hat{r}} < 3$ to hold for all $\hat{r} \in I_{d}$. For $d \geq 4$, we have $2 \cdot \floor*{\frac{n}{d}} - 3 \leq 2 \cdot \floor*{\frac{n}{4}} - 3 \leq 2 \cdot \frac{n}{4} - 3 = \frac{n}{2} - \frac{6}{2} = \frac{n - 6}{2} < \frac{n - 4}{2}$. Thus the $\frac{n - 4}{2}$ moves of $B$ suffice to remove enough numbers with residue $0$ and residue $\frac{d}{2}$.

For $d = 2$, there are equally many numbers with residue $0$ as with residue $1$. Therefore, the numbers on the board can be paired up as even and odd numbers, such that each number belongs to exactly one pair. Then  $B$ can win by the following strategy: For every even number that $A$ removes, $B$ removes an odd number, and vice versa. This way, only pairs consisting of one even and one odd number are removed. Since initially only such pairs exist, one such pair will remain at the end, so the sum of the final two numbers is odd and thus not divisible by $d = 2$.
\end{proof}

So far, we have only considered strategies that guarantee a win for one player, regardless of the opponent’s moves. However, it is also possible for $B$ to force a win in the game $\textit{Z}(n, d)$ for some odd $n \geq 5$ with $d \in \{\frac{n - 1}{2}, \frac{n + 1}{2}, n, n + 1, n + 2 \}$ and $d \neq 2$, if $A$ fails to play optimally. 
\vspace{\lineskip}

The restriction $d \neq 2$ is only necessary for $n = 5$ and $d = \frac{n - 1}{2} = \frac{5 - 1}{2} = 2$, because for $d = 2$, $A$ can secure a win even if $A$ does not play optimally.

\begin{theorem} \label{3.2.6}
Let $n \geq 5$ be odd. If $A$ fails to play optimally, $B$ can guarantee a win in the game $\textit{Z}(n, d)$ for $d \in \big \{\frac{n - 1}{2}, \frac{n + 1}{2}, n, n + 1, n + 2 \big \}$ with $d \neq 2$.
\end{theorem}

\begin{proof}
$B$ pursues the following strategy: For every $i \in \{0, \dots, d - 1\} \setminus I_{d}$ with $1 \leq a_{i} \leq a_{d - i}$, $B$ removes numbers with residue $i$ modulo $d$ until $a_{i} = 0$. For every $\hat{r} \in I_{d}$ with $a_{\hat{r}} \geq 2$, $B$ removes numbers with residue $\hat{r}$ until $a_{\hat{r}} < 2$. It is evident that this is the most efficient way for $B$ to create superfluous numbers.

It remains to be shown that $B$ can make all numbers superfluous before $A$ makes the final move of the game. Since $B$ has a total of $\frac{n - 3}{2}$ moves, it suffices to show that $B$ needs at most $\frac{n - 3}{2}$ moves to make all numbers superfluous.

Since $A$ is assumed not to play optimally, there is a game state in which $A$ is to move and could create an $A$-situation, but fails to do so. That is, $A$ removes a number whose residue $r$ modulo $d$ already satisfies the condition for an $A$-situation.
In the case $r \in I_{d}$, this implies that $a_{r}$ decreases by $1$ and becomes odd; in the case $r \notin I_{d}$, we then have $a_{r} = a_{d - r} - 1$. In either case, due to $A$’s suboptimal move, $B$ needs at least one fewer move to make all remaining numbers superfluous. 

We consider two cases for $d$.
\begin{enumerate}
    \item Let $d = \frac{n \pm 1}{2}$. If $d$ is even, then $d \geq 4$, $\lvert I_{d} \rvert = 2$, and thus $\lvert \{0, \dots, d - 1\} \setminus I_{d} \rvert = d - 2 \geq 4 - 2 = 2 > 0$. Initially, there are $\frac{\lvert \{0, \dots, d - 1\} \setminus I_{d} \rvert}{2} = \frac{d - 2}{2}$ pairs of residues $r \notin I_{d}$ modulo $d$, each with their additive inverse $d - r$ modulo $d$.
    
    First, let $d = \frac{n - 1}{2}$. As shown in \hyperref[1. of 3.1.4]{\ref{1. of 3.1.4}.} in the proof of \hyperref[3.1.4]{Theorem \ref*{3.1.4}}, we have $a_{r} = 2 < 3 = a_{1}$ for all $r \in \{0, 2, \dots, d - 1\}$. Therefore, for each of the $\frac{d - 2}{2}$ pairs of residues $r \notin I_{d}$ with their respective additive inverse $d - r$, $B$ must remove two numbers to make the corresponding numbers superfluous. Thus $B$ must remove $\frac{d - 2}{2} \cdot 2 = d - 2$ numbers. 
    Additionally, $B$ must remove one number with residue $0$ and one with residue $\frac{d}{2}$ to make the numbers with residue $\hat{r} \in I_{d}$ superfluous. Hence $B$ would need at most $(d - 2) + 2 = d = \frac{n - 1}{2}$ moves to make all numbers superfluous.
    \vspace{\lineskip}
    
    If $d = \frac{n + 1}{2}$, then $n \bmod \frac{n + 1}{2} = \frac{n - 1}{2} = d - 1$, so that $a_{0} = \floor*{\frac{n}{\frac{n + 1}{2}}} = \floor*{\frac{2n}{n + 1}} = 1$ and $a_{r} = 1 + 1 = 2$ for all $r \in \{1, 2, \dots, d - 1\}$. Again, there are $\frac{d - 2}{2}$ pairs of residues $r \notin I_{d}$ together with their additive inverse $d - r$, for which $B$ has to remove two numbers each to make the corresponding numbers superfluous.   
    Since $a_{0} = 1$ and $a_{\frac{d}{2}} = 2$, $B$ only needs to remove one number with residue $\frac{d}{2}$ to make all numbers with a residue $\hat{r} \in I_{d}$ superfluous.   
    Thus $B$ would again need at most $\frac{d - 2}{2} \cdot 2 + 1 = d - 1 = \frac{n + 1}{2} - 1 = \frac{n - 1}{2}$ moves to make all numbers superfluous. 
    \vspace{\lineskip}

    If $d$ is odd, then $d \geq 3$, $I_{d} = \{0\}$, and thus $\lvert \{0, \dots, d - 1\} \setminus I_{d} \rvert = d - 1 \geq 3 - 1 = 2 > 0$.  
    Regardless of whether $d = \frac{n - 1}{2}$ or $d = \frac{n + 1}{2}$, $B$ needs -- analogously to above -- at most $\frac{\lvert \{0, \dots, d - 1\} \setminus I_{d} \rvert}{2} \cdot 2 = \lvert \{0, \dots, d - 1\} \setminus I_{d} \rvert = d - 1$ moves to make all numbers with residue $r \notin I_{d}$ superfluous.   
    If $d = \frac{n - 1}{2}$, there are two numbers with residue $0$, of which $B$ only needs to remove one.  Thus $B$ would again need at most $(d - 1) + 1 = d = \frac{n - 1}{2}$ moves.   
    If $d = \frac{n + 1}{2}$, then $d$, having residue $0$, is the only number with a residue lying in $I_{d}$, so it is already superfluous from the very beginning.   
    Hence $B$ would again need only $d - 1 = \frac{n + 1}{2} - 1 = \frac{n - 1}{2}$ moves to make all numbers superfluous. 
    \vspace{\lineskip}

    In either case, $B$ would thus need at most $\frac{n - 1}{2}$ moves to make all numbers superfluous. However, since $A$ does not play optimally, $B$ needs at least one fewer move, that is, at most $\frac{n - 1}{2} - 1 = \frac{n - 3}{2}$ moves to make all numbers superfluous. Therefore, $B$ has enough moves to secure a win if $A$ does not play optimally.
    
    \item Let $d = n$ or $d = n + 2$. Then $d$ is odd, $I_{d} = \{0\}$, $\lvert \{0, \dots, d - 1\} \setminus I_{d} \rvert = d - 1$, and $a_{r} \leq 1$ for all $r \in \{0, \dots, d - 1\}$.  
    If $d = n$, then $n$ is superfluous, and $a_{r} = 1 = a_{d - r}$ for all $r \in \{1, \dots, n - 1\}$, so there are $\frac{d - 1}{2} = \frac{n - 1}{2}$ pairs of residues whose corresponding numbers on the board are not superfluous.  
    
    If $d = n + 2$, then $1$ is superfluous since $a_{1} = 1$ and $a_{d - 1} = a_{(n + 2) - 1} = a_{n + 1} = 0$. Also, $a_{0} = 0$ and $a_{r} = 1 = a_{d - r}$ for all $r \in \{2, \dots, n\}$, so again there are $\frac{\lvert \{2, \dots, n\} \rvert}{2} = \frac{n - 1}{2}$ pairs of residues whose corresponding numbers on the board are not superfluous.
    \vspace{\lineskip}

    For each of the $\frac{n - 1}{2}$ pairs, $B$ needs exactly one move to make the corresponding numbers superfluous. Thus $B$ would need at most $\frac{n - 1}{2}$ moves until all numbers are superfluous. However, since $A$ does not play optimally, $B$ again needs at most $\frac{n - 1}{2} - 1 = \frac{n - 3}{2}$ moves. As a result, $B$ can force a win if $A$ fails to play optimally.
\end{enumerate}
The statement for $d = n + 1$ follows immediately from the statement for $d = \frac{n + 1}{2}$, because if $B$ can make sure that the sum of the final two numbers is not divisible by $\frac{n + 1}{2}$, it is, in particular, not divisible by $n + 1$.
\end{proof}

The statement of \hyperref[3.2.6]{Theorem~\ref*{3.2.6}} forms the theoretical foundation of my website for the game \textit{Zahlenschlacht}. The program on the website, against which one plays as Player $A$, follows a simple algorithm: in each turn, it looks for numbers that are not superfluous, removes them in a way that creates superfluous numbers as efficiently as possible, and repeats this process until all numbers are superfluous -- after that, it removes numbers at random. For the game variants implemented on the website, the following holds: if Player $A$ plays optimally, $A$ can force a win. However, if $A$ fails to play optimally, $B$ can force a win. Hence the program offers a practical demonstration of this property. The program can only be defeated by optimal play.  

In total, there are $57$ game variants available on the website. The values of $n$ that can be selected are the odd numbers satisfying $5 \leq n \leq 25$. Moreover, the website also includes variants in which the corresponding value of $d$ does not belong to the set $\left\{\frac{n - 1}{2}, \frac{n + 1}{2}, n, n + 1, n + 2\right\}$. 
\vspace{\lineskip}

Furthermore, there is a mode for playing against another human player who is physically present rather than playing against the program. The website allows to change between both modes.
\vspace{\lineskip}

My website for the game Zahlenschlacht can be found in \cite{num4}.

%% file: Section_4-Conclusion/Conclusion.tex
\section{Conclusion} \label{Section 4}

Inspired by a problem from the German National Mathematics Competition concerning a mathematical game, we developed and analyzed a comprehensive generalization of that game, which I have named \textit{Zahlenschlacht}, and identified optimal strategies for it. The simplicity of the gameplay is remarkable, as it stands in contrast to the surprisingly complex mathematical theory of the game.
\vspace{\lineskip}

For the game \textit{Z}$(n, d)$ with an odd $n \geq 5$, we showed that Player $A$ can guarantee a win for all $d \in \{2, 3, \frac{n - 1}{2}, \frac{n + 1}{2}, n, n + 1, n + 2\}$, and additionally for $d \in \{4, 5, 6\}$ when $n \geq 11$. Player $B$ can guarantee a win for all $d$ with $\frac{n + 3}{2} \leq d \leq n - 1$ or $d \geq n + 3$. Furthermore, we showed that, for even $n$, there is not a single variant where $A$ can secure a win as $B$ can always win here. Additionally, we found formulas that allow to generate infinitely many game variants for fixed $d$, for which $A$ can force a win. Moreover, I developed a website where selected variants of Zahlenschlacht can be played against my program or a human opponent, allowing players to apply the strategies justified in this paper.
\vspace{\lineskip}

Lastly, we state some directions for future research.
\begin{enumerate}[label=(\roman*)]
    \item Are there discernible patterns for the values of $d$ with $d \in \{7, \dots, \frac{n - 3}{2}\}$ for odd $n$, which determine whether $A$ or $B$ can secure a win?
    \item Let $n$ be even and let the roles of Player $A$ and Player $B$ be reversed: $A$ still begins the game, but now wins if the sum of the final two numbers is not divisible by $d$, while $B$ wins if the sum is divisible by $d$.
    \vspace{\lineskip}
    
    Are there interesting strategies that guarantee a win for one player?
\end{enumerate}